\documentclass[a4paper]{article}
\usepackage{xcolor}
\usepackage{amsmath,amssymb,amsthm}
\usepackage{graphicx}
\usepackage[colorlinks=true, allcolors=blue]{hyperref}
\usepackage{subfigure}
\usepackage{circuitikz}

\newcommand{\R}{\mathbb{R}}
\DeclareMathOperator{\rk}{rk}
\DeclareMathOperator{\re}{Re}

\DeclareMathOperator{\diag}{diag}

{\left(\begin{smallmatrix}}
{\end{smallmatrix}\right)}

\newenvironment{smallbmatrix}
{\left[\begin{smallmatrix}}
{\end{smallmatrix}\right]}

\theoremstyle{plain}
\newtheorem{defi}{Definition}[section]
\newtheorem{theorem}[defi]{Theorem}
\newtheorem{prop}[defi]{Proposition}
\newtheorem{lemma}[defi]{Lemma}

\newtheorem{remark}[defi]{Remark}

\theoremstyle{definition}

\title{Port-Hamiltonian Modeling and Control of Electric Vehicle Charging Stations}

\author{Hannes~Gernandt, \thanks{H. Gernandt is with the Department
of Mathematics, Bergische Universität Wuppertal, 42119 Wuppertal, Germany e-mail: gernandt@uni-wuppertal.de and with the Fraunhofer Research Insitution for Energy Infrastructures and Geothermal Systems IEG, 03046 Cottbus, Germany.}
\and 
        Bernardo~Severino$^\ddagger$, \and Xinyi Zhang$^\ddagger$, \and  Volker~Mehrmann \thanks{V. Mehrmann is with the Department
of Mathematics, Technische Universität Berlin, 10623 Berlin,  Germany e-mail: mehrmann@math.tu-berlin.de.}  \and Kai~Strunz
\thanks{B.~Severino, K.~Strunz, X.~Zhang are with the SENSE Laboratory, Department of Electrical Engineering and Computer Sciences, Technische Universität Berlin, 10587 Berlin,
Germany (e-mail: kai.strunz@tu-berlin.de)}
}

\begin{document}

\maketitle

\begin{abstract}
Electric vehicles (EV) are an important part of future sustainable transportation. The increasing integration of EV charging stations (EVCSs) in the existing power grids require new scaleable control algorithms that maintain the stability and resilience of the grid. Here, we present such a control approach using an averaged port-Hamiltonian model. In this approach, the underlying switching behavior of the power converters is approximated by an averaged non-linear system. The averaged models are used to derive various types of stabilizing controllers, including the typically used PI controllers. The pH modeling is showcased by means of a generic setup of an EVCS, where the battery of the vehicle is connected to an AC grid via power lines, converters, and filters. Finally, the control design methods are compared for the averaged pH system and validated using a simulation model of the switched charging station.
\end{abstract}
 \vspace{0.5cm}\\


\section{Introduction}

The widespread adoption of electric vehicles (EVs) depends heavily on the extensive development of charging infrastructure and associated technologies. There is an expectation of a significant increase in the number of EV charging stations (EVCS) in the near future. Deploying an EVCS requires AC/DC converters for power conversion and DC/DC converters to regulate the charging process of EVs \cite{SAE3068}. Even though an EV's charging power may be adapted based on its state of charge (SOC), for transient stability analysis, the EV, with its DC/DC converter, is modeled as a constant power load (CPL), since the SOC remains relatively stable over the short periods considered \cite{1658410}. CPLs exhibit nonlinear behavior and a negative impedance $V$-$I$ characteristic, posing challenges to voltage stability \cite{8767007}. Therefore, managing the potential deterioration of voltage stability due to CPLs is vital for ensuring the robust operation of an EVCS.

A standard control approach to ensure robust operation of power systems in general are proportional-integral (PI) or cascaded PI controllers~\cite{Blaa06,YazI10}, which were particularly applied to EVCS in \cite{AranStru12} to achieve a robust control of the battery voltage. For a~particular model of the charging station, these controllers provide good results, and an appropriate choice of control gains of these controllers can be obtained from the linearized system dynamics around the considered operation point. However, due to the switching dynamics, the performance of the PI controllers may deviate from expectations when subjected to large disturbances. Another limitation of PI controllers arises when an additional component is introduced into the existing system. In such a~case, the derivation of new models to compute the required transfer functions for the control design process is needed.

To deal with these drawbacks, in recent years passivity-based and port-Hamiltonian-based nonlinear control has been considered. Those approaches are based on a physics-inspired  modeling of the system which explicitly takes into account the system's energy. In particular, port-Hamiltonian (pH) modeling of dynamical systems has recently received increased attention because of its robustness, its modeling flexibility, its ability to model multi-physical systems, and its advantageous properties for stability- and passivity-preserving model reduction, numerical simulation, and optimal control methods \cite{MehU22,SchJ14}. Furthermore, the pH modeling paradigm comprises the decomposition of the overall system model into several pH subsystems, allowing for a systematic built-up of the whole system. Due to these properties, pH control is a suitable choice for conducting real-time dynamical simulation and control of pH models for power systems.

Previously, pH modeling has been applied to power networks, see also \cite{SchZOSSR16} for a survey on the modeling of microgrids. In particular, various pH formulations of the synchronous generator are given in \cite{FiaZOSS12,FiaZOSS13,Ste18,SchS16} ranging from a scalar second-order swing equation to more accurate higher dimensional models. Furthermore, there exist pH models for the phase-locked loop (PLL), which is a standard component of power converters used to synchronize the output voltage of the point of common coupling (PCC) with the grid voltage \cite{Brav18}. 

When studying pH models, only few authors consider the switched system dynamics \cite{EscSO99}, and in the context of pH control design, mostly the averaged model instead of the switched system is used. PH control is applied to islanded AC microgrids in \cite{StreNMHF21} and more general types of grids \cite{ZhoS22}. Furthermore, \cite{MonMOS19} provides a \emph{shifted} pH formulation of power networks that contain constant power loads to estimate the regions of attraction. Additionally, averaged pH systems were recently used for observer design of DC power systems that contain DC/DC converters \cite{MegPL13}.

As an alternative to the pH control approach, shifted and Krasovskii passivity has been exploited for control of power systems in \cite{KosaCSP21}, but neither PI controllers, nor coupling to an AC grid, nor switched models of the power converters were considered to validate the passivity-based controller.

The main contributions of the present work are the following:
\begin{itemize}
    \item description of a general class of nonlinear systems of \emph{averaged pH systems} that can be used to model switched power systems;
    \item derivation of simple non-linear proportional (P) and proportional-integral (PI) control laws based on the averaged pH structure whose two control gains can be tuned easily; 
    \item showcasing the pH modeling approach for a~generic model of an EVCS leading to an averaged pH system that can easily be adapted to model variations;  
    \item verification of the proposed pH controllers that are designed for the averaged model;
    \item validation and comparison with a standard cascaded PI controller in a simulation of the original switching system dynamics.
\end{itemize}

The proposed generic pH model of an EVCS is based on the dynamic model which was derived in \cite{AranStru12} and consists of five subsystems, all described via a~pH model. These are (i) a~load model, (ii) a~DC/DC converter model, (iii) a short line model, (iv) an AC/DC converter model, as well as (v) PI-controller models. After putting the generic model of the EVCS in the pH framework, it can be easily integrated into existing pH power grid models or combined with several pH charging station models.

The paper is structured as follows: in Section~\ref{sec:ph}, we review pH modeling with a particular focus on the modeling of switching power systems. In Section~\ref{sec:control}, we analyze the steady states of the systems and propose pH control designs leading to asymptotically stable pH systems. In Section~\ref{sec:charging}, we derive a generic averaged pH model of an EVCS. Finally, we apply and compare the proposed pH controllers to the charging station in Section~\ref{sec:comp}. Proofs of the main results are presented in the Appendix \ref{sec:app_control} and \ref{sec:app_hurwitz}.

\section{PH state space models}
\label{sec:ph}
We derive a  pH state space model of a charging station. Because the pH modeling of the other components is relatively standard, we focus on the switching power systems arising in the converters. The complete pH model allows us to conduct a stability and passivity  analysis and to study control design for these systems. 

\subsection{PH models for  subsystems}\label{sec:subsystems}
To generate a complete model for the charging station, we first model each subsystem as a pH system using the model description introduced in \cite{Sch17} 
\begin{align}
\label{pH_DAE}
\begin{split}
\dot x&=~~\,(J-R)(x)\nabla \mathcal{H}(x)+(B-P)(x)u,\\
y&=(B+P)^\top(x)\nabla\mathcal{H}(x)+(S-N)(x)u,
\end{split}
\end{align}
for some matrix functions $J,R,S,N:\mathbb R^n\rightarrow\mathbb R^{n\times n}$, $B,P:\mathbb R^n\rightarrow\mathbb R^{n\times m}$, for a time interval $ \mathbb I$, the system state $x:\mathbb I\to \mathbb R^n$, the input $u:\mathbb I\to \mathbb R^m$, and the output $y:\mathbb I\to \mathbb R^m$. Here $\mathcal{H}$ is the \emph{Hamiltonian}, that describes the internal energy of the system and is  assumed to be a quadratic function of the form $\mathcal{H}(x)=\tfrac12 x^\top Hx$ for some symmetric positive definite  $H\in\mathbb{R}^{n\times n}$.  Furthermore, we assume for all states $x\in\mathbb{R}^n$ that 
\begin{align*}
&J(x)=-J(x)^\top,\ N(x)=-N(x)^\top,\\  &\begin{bmatrix}
\!\!R(x)&\!\!\!\!P(x)\\~P(x)^\top&\!\!\!\!S(x)
\end{bmatrix}\geq 0,
\end{align*}
where for a symmetric matrix function $M=M^\top$, $M\geq 0$  ($M>0$) denotes that the matrix is positive semidefinite (definite). 
Generalizations of pH systems to the descriptor case and to systems with non-quadratic Hamiltonians were given in \cite{MehM19}, see also \cite{MehU22} for an extensive survey.

An important property of PH systems is that they are \emph{passive}, meaning that they fulfill the \emph{power balance equation} and the \emph{dissipation inequality}
\begin{align}
&~~~~\dot{\mathcal{H}}(x)\nonumber \\&=u^\top y-\begin{bmatrix}
\nabla\mathcal{H}(x)\\u
\end{bmatrix}^\top\!\!\begin{bmatrix}
R(x)&\!\!\!\!P(x)\\P(x)^\top&\!\!\!\!S(x)
\end{bmatrix}\!\!\begin{bmatrix}
\nabla\mathcal{H}(x)\\u
\end{bmatrix}\label{ineq:dissi}\\&\leq u^\top y, \nonumber
\end{align}
where the time derivative $\dot{\mathcal{H}}$ is usually interpreted as the \emph{stored power}, $u^\top y$ is interpreted as the \emph{supplied power} and the remaining summand is viewed as the \emph{dissipated power}. 

If $u=y=0$ and if for a stationary point $\overline{x}$ the matrix $(J-R)(\overline{x})$ is \emph{Hurwitz}, i.e., all eigenvalues are in the open left half-plane, then $\mathcal H$ is a Lyapunov function implying that the system is \emph{locally exponentially stable}, see \cite[Theorem 5.35]{LogR14}. Note, however, that if $R(\overline{x})$ is not positive definite, then we can in general not use $\mathcal{H}$ as a Lyapunov function, but have to first transform the system, see \cite{AchAM21}.

The components of the charging station, such as converters and lines, have in-going current/voltage $i_i,v_i:\mathbb I \to \mathbb R^{k_i}$, as well as out-going current/voltage $i_o,v_o:\mathbb I \to \mathbb R^{k_o}$. This special structure of the components leads to the following pH subsystems
\begin{align}
\label{special_ph}
\begin{split}
\dot x&=(J(x)-R(x))\nabla \mathcal{H}(x)+\begin{bmatrix}
B_{i}&B_{o}
\end{bmatrix}\begin{bmatrix}u_{1}\\ u_{2}\end{bmatrix},\\
y&=\begin{bmatrix}y_{1}\\y_{2}\end{bmatrix}=\begin{bmatrix}
B_{i}^\top\\B_{o}^\top
\end{bmatrix}\nabla\mathcal{H}(x)
\end{split}
\end{align}
for some $B_i\in\R^{n\times k_i}$, $B_o\in\R^{n\times k_o}$ and with 
\begin{align}
\label{eq:yu}
y^\top u=y_{1}u_{1}+y_{2}u_{2}=i_{i}^\top v_{i}-i_{o}^\top v_{o}.
\end{align}
Using this particular input-output structure of the pH system, the dissipation inequality \eqref{ineq:dissi} implies that the out-going power of the components is less than or equal to the in-going power. However, \eqref{eq:yu} provides some freedom of choice in the input and output variables. Typical choices are 
\begin{itemize}
    \item[(i)] \emph{voltage input} systems, where $(u_{1},u_{2})=(v_{i},v_{o})$ and $(y_{1},y_{2})=(i_{i},-i_{o})$;
    \item[(ii)] \emph{current input} systems, where $(u_{1},u_{2})=(i_{i},-i_{o})$ and $(y_{1},y_{2})=(v_{i},v_{o})$;
    \item[(iii)] \emph{current/voltage input} systems, where $(u_{1},u_{2})=(v_{i},-i_{o})$ and $(y_{1},y_{2})=(i_{i},v_{o})$, or vice versa.
\end{itemize}
\subsection{Interconnection of pH systems}\label{sec:interconn}
An important property of  pH representations is that for two given pH systems, their \emph{power conserving interconnection} will again be pH and, in particular passive, see e.g.\ \cite[Section 2.2.4]{MehM19} and \cite[Chapter 6]{Sch17}. In the context of power systems, this interconnection is described by Kirchhoff's current and voltage laws. For interconnection, we consider two subsystems of the form 
\begin{align}
\label{special_ph_2}
\begin{split}
\dot x_j&=(J_j-R_j)(x_j)\nabla \mathcal{H}_j(x_j)+\begin{bmatrix}
B_{i,j}&\!\!\!B_{o,j}
\end{bmatrix}\!\begin{bmatrix}u_{i,j}\\ u_{o,j}\end{bmatrix},\\
y&=\begin{bmatrix}y_{i,j}\\y_{o,j}\end{bmatrix}=\begin{smallbmatrix}
B_{i,j}^\top\\B_{o,j}^\top
\end{smallbmatrix}\nabla\mathcal{H}_j(x_j),\quad j=1,2.
\end{split}
\end{align}

If we assume that the direction of the current flow and the sign of the voltages is the same for both subsystems, the natural coupling conditions based on Kirchhoff's laws are
\begin{align}
\label{VI_coupl}
i_{o,1}=i_{i,2},\quad v_{o,1}=v_{i,2}.
\end{align}
In the following, we will present a typical coupling of two power systems which  both have current and voltage inputs. Let the subsystems be given by 
\begin{align*}
\dot x_j&=(J_j-R_j)(x_j)\nabla \mathcal{H}_j(x_j)+\begin{bmatrix}
B_{i,j}&\!\!\!\!B_{o,j}
\end{bmatrix}\!\begin{bmatrix}v_{i,j}\\ -i_{o,j}\end{bmatrix}\!,\\
y_j&=\begin{bmatrix}i_{i,j}\\v_{o,j}\end{bmatrix}=\begin{bmatrix}
B_{i,j}^\top\\[1ex]B_{o,j}^\top
\end{bmatrix}\nabla\mathcal{H}_j(x_j),\quad j=1,2.
\end{align*}
Then, the interconnected system based on \eqref{VI_coupl} is given by 
\begin{align*}
\begin{bmatrix}
\dot x_1\\ \dot x_2
\end{bmatrix}&\!=\!\begin{bmatrix}
(J_1-R_1)(x_1)&\!\!\! - B_{o,1}
B_{i,2}^\top\\ 
B_{i,2}B_{o,1}^\top&\!\!\!(J_2-R_2)(x_2)\end{bmatrix}\!\begin{bmatrix}
\nabla\mathcal{H}_1(x_1)\\\nabla\mathcal{H}_2(x_2)
\end{bmatrix} \\&~~~+\begin{bmatrix}B_{i,1}&\!\!\!0\\0&\!\!\!B_{o,2}\end{bmatrix}\begin{bmatrix}
v_{i,1}\\-i_{o,2}
\end{bmatrix},\\[1ex]
y&=\begin{bmatrix}
i_{i,1}\\v_{o,2}
\end{bmatrix}=\begin{bmatrix}B_{i,1}^\top&\!\!\!0\\0&\!\!\!B_{o,2}^\top\end{bmatrix}\begin{bmatrix}
\nabla\mathcal{H}_1(x_1)\\\nabla\mathcal{H}_2(x_2)
\end{bmatrix}.
\end{align*}
Couplings for the other system types (i) or (ii) can be obtained analogously. More general couplings of more than two subsystems can be described using so called \emph{Dirac structures} \cite[Section 6.6]{Sch17}. 

\subsection{Averaged pH models}\label{sec:averaged}
The considered power electronic systems in a~charging station have fast switching components which are parts of the DC/DC and the AC/DC converters. For the purpose of control system design and analysis, we approximate the switching behavior and use state space averaged pH models. 
For a detailed modeling of switching pH systems, see \cite{EscSO99}. We consider a~periodic switching between two given pH systems of the form
\[
\dot x=(J_j-R_j)(x)\nabla\mathcal{H}(x)+Bu_j,\quad j=1,2,
\]
where we assume that $Bu=B_1u_1=B_2u_2$, which means that the inputs of the switched subsystems are the same. Although the averaged pH system can still be formulated without this assumption, it is used to shorten the presentation.

The \emph{duty ratio} $d:\mathbb I \to [0,1]$ describes the average by which the two systems identified by $j=1$ or $j=2$ are weighted, respectively. In particular, $d=0$ means that we keep the model $j=1$ for all $t\in \mathbb I$ and do not switch. On the other hand, $d=1$ means that we keep the model $j=2$, $t\in \mathbb I$. Then, the \emph{state space averaged model} which approximates the switching system is given by 
\[
\dot x=((1-d)(J_1-R_1)+d(J_2-R_2))(x)\nabla\mathcal{H}(x)+Bu.
\]
Note that in the literature, see e.g. \cite{JudCMCG17}, sometimes the components are also modeled by equations using the \emph{modulation index} $m\in[-1,1]$, which is related to the duty ratio via $d=\tfrac{m+1}{2}$. 

While standard pH systems are well understood in terms of stability and passivity, see e.g. \cite{MehU22,SchJ14}, for averaged pH systems this is new and will be discussed in the next section.

\section{Stability analysis and control design for averaged pH systems} \label{sec:control}
In this section, we determine the steady states of averaged pH systems described in Section \ref{sec:averaged}, and study the stability of systems which include the duty ratios as control variables. 

To this end, we assume that for $k$ switched subsystems with parameters $\theta_j\in[0,1]$, $j=1,\ldots, k$ the averaged pH system has the form 
\begin{align}
\label{def_sys}
\begin{split}
\dot x&=D(\theta)\nabla\mathcal{H}(x)+Bu,\\
y&=B^\top \nabla\mathcal{H}(x),
\end{split}
\end{align}
where 
\begin{align*}
D(\theta)&:=J-R+\sum_{j=1}^k\theta_j(J_j-R_j),\\ 
D&:=J-R,\ D_j:=J_j-R_j,\ j=1,\ldots,k, \\ J&\phantom{:}=-J^\top,\ R=R^\top\geq 0\\
J_j&\phantom{:}=-J_j^\top,\ R_j=R_j^\top\geq 0,\ j=1,\ldots,k,
\end{align*}
are all matrices in $\mathbb R^{n\times n}$, $B\in\mathbb R^{n\times m}$, and the Hamiltonian is $\mathcal{H}(x)=\tfrac12x^\top Hx$ for some positive definite $H\in\mathbb R^{n\times n}$. 
Note that we can replace duty ratios by modulation indexes and vice versa, the resulting additional terms can also be combined in the matrices $J$ and $R$.

In the following subsection we study the equilibrium points and also the stabilization by controlling the duty ratios. The control design is based on extending the system by introducing $\theta$ as an additional state variable.

\subsection{Extended system and steady states}\label{sec:extsys}
System \eqref{def_sys} depends on the vector of duty ratios $\theta=[\theta_1,\ldots,\theta_k]^\top\in[0,1]^k$. In this section, we define an extended system which contains these additional parameters as a state variable and derive a characterization when this system is locally asymptotically stable. For this, the system equations \eqref{def_sys} are rewritten as 
\begin{align*}
\dot x&=\begin{bmatrix}D&\mathcal{D}(x)\end{bmatrix}\begin{bmatrix}\nabla \mathcal{H}(x)\\\theta\end{bmatrix}+Bu.
\end{align*}
with
\[
\mathcal{D}(x):=\begin{bmatrix}D_1\nabla\mathcal{H}(x)&\ldots &D_k\nabla\mathcal{H}(x)\end{bmatrix}\in\R^{n\times k}.
\]
We consider two possible controller designs, a dynamic controller
\begin{equation}\label{dyncon}
\dot \theta=K_H(x)\nabla \mathcal{H}(x)-K_{\theta}\theta,
\end{equation}
and a stationary controller
\begin{equation}\label{statcon}
0=K_H(x)\nabla \mathcal{H}(x)-K_{\theta}\theta,
\end{equation}
with an invertible design parameter matrix $K_{\theta}\in\mathbb R^{k\times k}$, and a state-dependent matrix function $K_H:\R^n\rightarrow\mathbb R^{k\times n}$. 

The use of \eqref{dyncon} results in the extended pH system
\begin{align}
\label{eq:ext_ode_ws}
\begin{bmatrix}
\dot x\\ \dot \theta
\end{bmatrix}&=\begin{bmatrix}D&\mathcal{D}(x)\\K_H(x)&-K_{\theta}\end{bmatrix}\nabla \mathcal{H}^e(x,\theta) +\begin{bmatrix}B\\0\end{bmatrix}u,\\
 y&=\begin{bmatrix}B\\0\end{bmatrix}^\top \nabla\mathcal{H}^e(x,\theta), \nonumber
\end{align}
whereas the use of \eqref{statcon} leads to an extended differential-algebraic pH system 
\begin{align}
\label{eq:ext_dae_ws}
\begin{bmatrix} I_n & \!\!\! 0 \\ 0 & \!\!\! 0 \end{bmatrix} \!\!\begin{bmatrix}
\dot x\\ \dot \theta
\end{bmatrix}&=\begin{bmatrix}D&\mathcal{D}(x)\\K_H(x)&-K_{\theta}\end{bmatrix}\!\nabla \mathcal{H}^e(x,\theta) +\begin{bmatrix}B\\0\end{bmatrix}u,\\
 y&=\begin{bmatrix}B\\0\end{bmatrix}^\top \nabla\mathcal{H}^e(x,\theta). \nonumber
\end{align}
In both cases the extended Hamiltonian is
\[
\mathcal H^e(x,\theta)=\frac 12 \begin{bmatrix} x\\ \theta
\end{bmatrix}^\top \begin{bmatrix} H& 0 \\ 0 & I_k
\end{bmatrix}
\begin{bmatrix} x\\ \theta
\end{bmatrix},
\]
and the extended $J^e(x,\theta)$ and $R^e(x,\theta)$ are defined as the skew-symmetric and symmetric part of
\[
\begin{bmatrix}J-R&\mathcal{D}(x)\\K_H(x)&-K_\theta\end{bmatrix}
\]
for which the control design has to be chosen in such a way that the extended dissipation matrix satisfies $R^e(x,\theta)\geq 0$. One of many possible choices to achieve this is to choose $K_H(x)=-\mathcal D(x)^\top$ and $K_\theta>0$.

Note that a choice with $K_\theta$ invertible already guarantees that the differential-algebraic system has differentiation index one, see \cite{KunM06} which (at least theoretically) allows to remove the algebraic equations by solving for $\theta$.

We immediately obtain that by this extension, the steady state does not change.
 \begin{lemma}
Consider an averaged pH system given by \eqref{def_sys} with $u=\overline{u}$ and $\mathcal{H}(x)=\tfrac12x^\top Hx$ for some positive definite $H\in\R^{n\times n}$. Then the extended systems \eqref{eq:ext_dae_ws} and \eqref{eq:ext_ode_ws} have the same steady states $\overline{x}$ and $\overline{\theta}$ and these can be obtained from the equation
\begin{align}
\label{eq:steady}
0&=D(\overline{\theta})\nabla \mathcal{H}(\overline{x}) +B\overline{u}.
\end{align}
Furthermore, if $D(\overline{\theta})$ is invertible then  
\begin{align}
\label{eq:ss_inv}
\overline{x}=-H^{-1}D(\overline{\theta})^{-1}B\overline{u}.
\end{align}
\end{lemma}

\subsection{Control design for shifted averaged pH  systems}\label{sec:designshift}
In this section, we consider steady states of the extended system \eqref{eq:ext_ode_ws} and \eqref{eq:ext_dae_ws} for a given fixed control input $u=\overline{u}$. This results in a non-zero equilibrium state $\overline{x}$. Since the Hamiltonian of the system is often strictly convex with global minimum at $x=0$ which is different from the equilibrium $\overline{x}$, it is common to consider the \emph{shifted Hamiltonian}, see e.g.\  \cite{MonMOS19,Sch17}, which is given by 
\begin{align}
\label{shifted_ham}
\begin{split}
\mathcal{S}(x)&:=\mathcal{H}(x)-\mathcal{H}(\overline{x})-(x-\overline{x})^\top\nabla\mathcal{H}(\overline{x})\\&\phantom{:}=\mathcal{H}(x-\overline{x}).
\end{split}
\end{align}
The expression for $\mathcal{S}$ is called \emph{Bregman divergence} in convex analysis \cite{BerHM99} and fulfills
\[
\nabla \mathcal{S}(x)=\nabla\mathcal{H}(x)-\nabla\mathcal{H}(\overline{x}).
\]
Using the shifted Hamiltonian \eqref{shifted_ham} and extending it appropriately as 
\[
\mathcal{S}^e(x,\theta)=\mathcal{H}^e(x-\overline x,\theta-\overline \theta),
\]
we have
\[
\nabla \mathcal{S}^e(x,\theta)=\nabla\mathcal{H}^e(x,\theta)-\nabla\mathcal{H}^e(\overline x,\overline \theta).
\]
We also extend the system equations. Subtracting \eqref{eq:steady} from the first equation in \eqref{eq:ext_ode_ws} leads to 
\begin{align}
\dot x&=D\nabla \mathcal{S}(x)+\mathcal{D}(x)\theta-\mathcal{D}(\overline{x})\overline{\theta}+B(u-\overline{u}) \nonumber\\
\nonumber &=D\nabla \mathcal{S}(x)+\sum_{j=1}^{k}\theta_jD_j\nabla\mathcal{S}(x)\\&~~~+\mathcal{D}(\overline{x})(\theta_j-\overline{\theta_j})+B(u-\overline{u})\nonumber\\
&=\begin{bmatrix}D(\theta)&\!\!\mathcal{D}(\overline{x})\end{bmatrix}\nabla \mathcal{S}^e(x,\theta)+B(u-\overline{u})\label{eq:minusbar}.
\end{align}
For the shifted system, we again consider two controllers with design parameters $K_{\theta}\in\mathbb R^{k\times k}$ and $K_H\in\R^{k\times n}$ given by
\begin{align}\label{dynshift}
\dot \theta=K_H\nabla \mathcal{S}(x)-K_{\theta}(\theta-\overline{\theta}),
\end{align}
and
\begin{align}\label{statshift}
 0=K_H\nabla \mathcal{S}(x)-K_{\theta}(\theta-\overline{\theta}),
 \end{align}
where $K_H\in\R^{k\times n}$ and $K_{\theta}\in\R^{k\times k}$ are assumed to be independent of the state variables.
The use of \eqref{dynshift} results in the system
\begin{align}
    \label{eq:ext_ode}
\begin{bmatrix}
\dot x\\ \dot \theta
\end{bmatrix}&=\begin{bmatrix}D(\theta)& \!\!\!\mathcal{D}(\overline{x})\\K_H&\!\!\!-K_{\theta}\end{bmatrix}\nabla \mathcal{S}^e(x,\theta) +\begin{bmatrix}B\\0\end{bmatrix}(u-\overline{u}),\\
 y&=\begin{bmatrix}B\\0\end{bmatrix}^\top \nabla\mathcal{S}^e(x,\theta), \nonumber
\end{align}
whereas the use of the second controller leads to a differential-algebraic pH system
\begin{align}
\nonumber
\begin{bmatrix} I_n & 0 \\ 0 & 0 \end{bmatrix} \begin{bmatrix}
\dot x\\ \dot \theta
\end{bmatrix}&=\begin{bmatrix}D(\theta)&\mathcal{D}(\overline{x})\\K_H&-K_{\theta}\end{bmatrix}\nabla \mathcal{S}^e(x,\theta) \\&~~~+\begin{bmatrix}B\\0\end{bmatrix}(u-\overline{u}),\label{eq:ext_dae}\\
 y&=\begin{bmatrix}B\\0\end{bmatrix}^\top \nabla\mathcal{S}^e(x,\theta) \nonumber,
\end{align}
and extended $J^e(x,\theta-\overline{\theta})$ and $R^e(x,\theta-\overline{\theta})$ defined as the skew-symmetric and symmetric part of
\[
\begin{bmatrix}D(\theta)&\mathcal{D}(x)\\K_H&-K_{\theta}\end{bmatrix}
\]
for which the control design has to be chosen in such a way that the extended dissipation matrix satisfies $R^e(x,\theta-\overline{\theta})\geq 0$, see Theorem~\ref{prop:control} below.

Many control designs of power systems use PI-controllers and cascades thereof. To allow for PI-controllers in our duty ratio control laws, we further extend \eqref{eq:ext_dae} by using an integral state variable
$\xi\in\mathbb{R}^k$ and a further extended shifted Hamiltonian 
\begin{align}
\label{eq:shifted_Ham}
\mathcal S^f(x,\theta,\xi)=\mathcal{H}^e(x-\overline x,\theta-\overline{\theta}) +\tfrac12 \xi^\top\xi,
\end{align}
as well as the further extended differential-algebraic system
\begin{align}
\nonumber
\begin{bmatrix}
\dot x\\ 0 \\\dot \xi
\end{bmatrix}&=\begin{bmatrix}D(\theta)&\mathcal{D}(\overline{x})&0\\K_H&-K_{\theta}&I_{k}\\K_{\xi}&0&0\end{bmatrix}\nabla \mathcal{S}^f(x,\theta,\xi)\label{eq:PI}\\&~~~+\begin{bmatrix}B\\0\\0\end{bmatrix}(u-\overline{u})\\
 y&=\begin{bmatrix}B^\top&0&0\end{bmatrix} \nabla\mathcal{S}^f(x,\theta,\xi). \nonumber
\end{align}
Note that we do not include $\overline{\xi}$ in the definition of the extended shifted Hamiltonian \eqref{eq:shifted_Ham}, since the second line in \eqref{eq:PI} implies for the steady state values $x=\overline x$, $\theta=\overline{\theta}$ that $\overline{\xi}=0$ holds.

In the following result, we show that with the so constructed controllers, we can obtain locally asymptotically stable pH systems. 
\begin{theorem}
\label{prop:control}
Given an averaged pH system of the form \eqref{def_sys} with equilibrium $(\overline{x},\overline{\theta},\overline{u})$ and let $K_{\theta}\in\R^{k\times k}$ be positive definite such that $D(\overline{\theta})-\mathcal{D}(\overline{x})K_{\theta}^{-1}\mathcal{D}(\overline{x})^\top$ is Hurwitz. Then the extended system \eqref{eq:ext_ode} is locally asymptotically stable and pH for the following parameters:
\begin{itemize}
    \item[(a)] $\theta=\overline{\theta}$ and $K_H=0$;
    \item[(b)] $K_H=-\mathcal{D}(\overline{x})^\top$ and  $K_{\theta}\in\R^{k\times k}$ positive definite. 
\end{itemize}
Moreover, the following systems are equivalent to a locally asymptotically stable and pH system:
\begin{itemize}
    \item[(c)] System \eqref{eq:ext_dae} with  $K_H=-\delta\mathcal{D}(\overline{x})^\top$ for all $K_{\theta}\in\R^{k\times k}$ positive definite and $\delta>0$;
    \item[(d)] System \eqref{eq:PI} with $K_H=-\delta\mathcal{D}(\overline{x})^\top$, $K_{\theta}\in\R^{k\times k}$ positive definite, $\delta>0$, and $K_{\xi}=-(\mathcal{D}(\overline{x})K_{\theta}^{-1})^\top$.
\end{itemize}
\end{theorem}
The proof of Theorem~\ref{prop:control} is presented in Section~\ref{sec:app_control} in the appendix.

To prove that the shifted closed loop system is locally asymptotically stable and pH, it is crucial to assume that there exists positive definite $K_{\theta}\in\R^{k\times k}$ such that $D(\overline{\theta})-\mathcal{D}(\overline{x})K_{\theta}^{-1}\mathcal{D}(\overline{x})^\top$ is Hurwitz. In the remainder of this section, we discuss several sufficient conditions for this  to hold.

One trivial observation is that if the equilibrium state $\overline{\theta}$ is already chosen in such a way that $D(\overline{\theta})$ is Hurwitz, then any positive definite $K_{\theta}\in\R^{k\times k}$ leads to a Hurwitz matrix $D(\overline{\theta})-\mathcal{D}(\overline{x})K_{\theta}^{-1}\mathcal{D}(\overline{x})^\top$. 

Another recent characterization of $D(\overline{\theta})-\mathcal{D}K_{\theta}^{-1}\mathcal{D}^\top$ being Hurwitz was obtained in \cite{AchAM21} and it is based on considering its skew-symmetric part $A_J:=J(\overline{\theta})$ and its symmetric part  $A_R:=-R(\overline{\theta})-\mathcal{D}K_{\theta}^{-1}\mathcal{D}^\top$. One of several equivalent characterizations given in \cite[Lemma 3]{AchAM21} is that the matrix  $[A_R,A_JA_R,\ldots,A_J^{n-1}A_R]$ has full rank. 

In the next result, it is shown that the existence of $K_{\theta}$ satisfying the assumptions of Theorem \ref{prop:control} is guaranteed if the matrix $[R(\overline{\theta}),\mathcal{D}(\overline{x})]$ has full rank.
\begin{prop}
\label{lem:Hurwitz}
Let $J,R\in\R^{n\times n}$ satisfy $J=-J^\top$ and $R\geq 0$. Furthermore, let  $\mathcal{D}\in\R^{n\times k}$ be such that $\rk[R,\mathcal{D}]=n$ and $K_{\theta}\in\R^{k\times k}$ positive definite. Then $A:=J-R-\mathcal{D}K_{\theta}^{-1}\mathcal{D}^\top$ is Hurwitz.
\end{prop}
The proof of Proposition~\ref{lem:Hurwitz} is presented in Section~\ref{sec:app_hurwitz} in the appendix.

Another approach to obtain a locally asymptotically stable pH system for a given shifted averaged pH system \eqref{eq:minusbar} with equilibrium $(\overline{x},\overline{\theta},\overline{u})$ is based on the \emph{Bass stabilization algorithm}, where one chooses a solution $\hat K_{\theta}=\hat K_{\theta}^\top>0$ to the Lyapunov equation
\begin{align}
\label{eq:lyap_bass}
(\alpha I_n +A)\hat K_{\theta} + \hat K_{\theta}( \alpha I_n + A^\top) \ = \ 2\mathcal{D}(\overline{x})\mathcal{D}(\overline{x})^\top,
\end{align}
where $A:=D(\overline{\theta})$ and $\alpha>0$ has to exceed the largest singular value of $A$. Then the matrix $A-\mathcal{D}(\overline{x})\mathcal{D}(\overline{x})^\top\hat K_{\theta}^{-1}$ is Hurwitz.

It is well known, see e.g. \cite[p.\ 72]{TreSH01}, that a positive definite solution $\hat K_{\theta}$ of \eqref{eq:lyap_bass} exists if $(A,\mathcal{D}(\overline{x}))$ is stabilizable, i.e.\ $\rk[\lambda I_n-A,\mathcal{D}(\overline{x})]=n$ for all $\lambda\in\mathbb{C}$ satisfying ${\rm Re\,} \lambda\geq 0$. By assumption, $\lambda I_n-A$ is invertible for all $\lambda\in\mathbb{C}$ satisfying ${\rm Re\,} \lambda >0$. Therefore, the stabilizability of  $(A,\mathcal{D}(\overline{x}))$ is equivalent to 
\[
\rk[i\omega-D(\overline{\theta}),\mathcal{D}(\overline{x})]=n
\quad \text{for all $\omega\in\R$}.
\]
The positive definite solution $\hat K_{\theta}$ of \eqref{eq:lyap_bass} can then be used in the following equations
\begin{align*}
\dot \theta&=-\mathcal{D}(\overline{x})^\top\hat K_{\theta}^{-1}\nabla \mathcal{S}(x)-(\theta-\overline{\theta}),\\
 0&=-\mathcal{D}(\overline{x})^\top\hat K_{\theta}^{-1}\nabla \mathcal{S}(x)-(\theta-\overline{\theta}),
\end{align*}
which can alternatively be employed to extend the shifted averaged pH system. This extended system is locally asymptotically stable, which can be viewed as an analogue to the results in  Theorem~\ref{prop:control}, but to formally write the system as a pH system one has to multiply \eqref{eq:minusbar} with $\hat K_{\theta}$ which leads to a pH descriptor system, see \cite{MehM19} for further reading.
The control design which was outlined above was not needed to control the electric vehicle charging station because, for our application, the matrix $D(\overline{\theta})$ turned out to be Hurwitz already.

\begin{remark}
\label{rem:cascade}{\rm 
The control of EVCS often involves cascaded PI controllers \cite{AranStru12}. These controllers can be put in the proposed pH framework by further extending the system \eqref{eq:PI} and using $\hat u=u-\overline{u}$ leading to
\begin{align*}
\begin{bmatrix}\dot x\\ 0 \\\dot \xi_1\\ \dot \xi_2
\end{bmatrix}\!\!=\!\!\begin{bmatrix}D(\theta) &\!\!\!\!\mathcal{D}(\overline{x})&\!\!\!\!0&\!\!\!\!0\\K_H&\!\!\!\!-K_d&\!\!\!\!K_{\xi,1}&\!\!\!\!K_{\xi,2}\\K_{\xi,3}&\!\!\!\!0&\!\!\!\!0&\!\!\!\!0\\K_{\xi,4}&\!\!\!\!0&\!\!\!\!K_{\xi,5}&\!\!\!\!0
\end{bmatrix}\!\!\!\begin{bmatrix}\nabla \mathcal{S}(x)\\\theta-\overline{\theta}\\\xi_1\\ \xi_2\end{bmatrix}\!\!+\!\!\begin{bmatrix}
B\\ 0\\0 \\ 0 \end{bmatrix}\hat u.
\end{align*}
 }
\end{remark}

\section{EVCS model as (averaged) pH system}
\label{sec:charging}
\begin{center}
\begin{figure*}[t]
\centering
\includegraphics[scale=0.9]{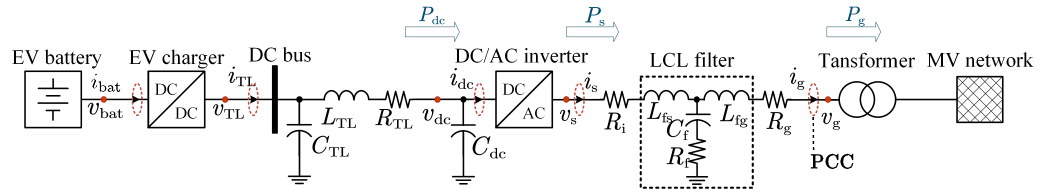}
\caption{Circuit diagram of the EV charging station. The bidirectional DC-DC converter and the EV battery model are shown in the Figures~\ref{bm} and \ref{bm2}.}
\label{cs}
\end{figure*}
\end{center}
Here, we derive a complete pH system model for the charging station without controllers, as shown in Figure~\ref{cs}.
The charging station system encompasses, from left to right, the battery of an EV, a DC/DC bidirectional converter, a DC bus, a DC line, a DC/AC voltage source converter (VSC), an LCL output filter, an AC transformer, and the connection to the medium voltage (MV) electrical grid. The boost converter is used to step up the DC output voltage of the EV battery to the DC bus voltage.
The central voltage source controller is used to interface the DC bus to the three-phase AC network. The transformer is utilized to elevate the low voltage in the charging station to the medium voltage of the MV network. Furthermore, an LCL filter is used to filter out the high frequency harmonics generated by the VSC. 

In the following, the pH model of the charging station is obtained by modeling each of the subsystems as pH systems and by interconnecting them. Here one usually considers variables which represent energy of the system, i.e.\ instead of voltages of capacitors and currents of inductors, we will use charges and fluxes. 

\subsection{DC side subsystem}
\label{sec:dc}

\subsubsection*{Battery model}
The battery model is shown in Figure~\ref{bm}, see \cite{7442180}. 
\begin{figure}[h]
\centering

\begin{circuitikz}[xscale=-1, american] 
    \draw (0,0) 
    to[open, v=$~~~~~~v_{\text{dc}}$]  (0,2) 
    to[short] (1,2)
    to[C=$C_{\text{TL}}$] (1,0) -- (0,0);
    \draw (1,2)
   to[closing switch=$S_1$] (3,2);
    
    \draw (3,2)
    to[L=$L$] (5,2);
    
    \draw (1,0)
    to (3,0);

    \draw (3.1,0)
    to[closing switch=$S_2$] (3.1,2);

    \draw (3,0)
    to[short] (5,0);
    
    \draw (5,2)
    to[C=$C$] (5,0);
    
    \draw (5,2)
    to[short] (6,2)
    to[open,v<=$v_{\text{bat}}$] (6,0)
    to[short] (5,0);
    
    \draw (0,0) to[short] (5,0);

    \draw (2.2,2)
    to[short, i=$i_L$] (3.3,2);
     \draw (5,2)
    to[short, i=$i_{bat}$] (6,2);
    \draw (0,2)
   to[short, i=$i_{dc}$] (1,2);

    \draw (0,0) to[short, -*] (0,0);
    \draw (0,2) to[short, -*] (0,2);
    \draw (6,0) to[short, -*] (6,0);
    \draw (6,2) to[short, -*] (6,2);
\end{circuitikz}
\caption{The DC-DC converter model used in Figure~\ref{cs}.} \label{bm}
\vspace{2ex}
\begin{circuitikz}[american voltages]
    \draw
    (0,0) to[cV, v<=$v_{ev}$] (0,2) 
    to[R, l=$R_{bat}$, i>=$i_{bat}$] (3,2) 
    to[short, -*] (3,2) 
    (3,0) to[short, *-] (3,0) 
    to (0,0); 
    \draw (3,2) to[open, v^>=$v_{bat}$, invert] (3,0); 
\end{circuitikz}
\caption{The EV battery model used in Figure~\ref{cs}.} \label{bm2}
\end{figure}
The open-circuit voltage is denoted by $v_{ev}$, and the voltage-current characteristic is 
modeled by a resistance $R_{bat}$.
The Kirchhoff laws imply that
\begin{align}
\label{eq:bat}
v_{bat}=v_{ev}+R_{bat}i_{bat},\quad 
-i_{ev}=i_{bat}.
\end{align}

\subsubsection*{DC/DC converter} 
For the charging station we use a bidirectional DC/DC converter as in \cite{6596508} that is located in Figure~\ref{cs} between the EV battery and the DC Bus. This subsystem connects to the battery via the port variables $(v_{bat},i_{bat})$ and to a DC bus via the port-variables  $(v_{\text{\tiny TL}},i_{\text{\tiny TL}})$.

As usual in pH modeling, we replace current $i_L$ by the flux $\phi_{L}=Li_{L}$ and voltage $v_{bat}$ by charge $Q_{bat}=Cv_{bat}$, $Q_{\text{\tiny TL}}=C_{\text{\tiny TL}}v_{\text{\tiny TL}}$. Then an averaged model of the converter is given by 
\begin{align}
\nonumber
  \dot x_{bdc}&=\begin{bmatrix}0&\!\!\!\!\!\!-1&\!\!\!d_{dc}\\ 1 &\!\!\!\!\!\!\phantom{-}0&\!\!\!0\\-d_{dc}&\!\!\!\!\!\!\phantom{-}0&\!\!\!0\end{bmatrix}\!\!H_{bdc}x_{bdc}\!+\!\begin{bmatrix}0&\!\!0\\1&\!\!0\\0&\!\!1\end{bmatrix}\begin{bmatrix}
 -i_{bat}\\i_{\text{\tiny TL}}
\end{bmatrix}\!,\\
y&=\begin{bmatrix}0&1&0\\0&0&1\end{bmatrix}H_{bdc}x_{bdc}=\begin{bmatrix}
v_{bat}\\v_{\text{\tiny TL}}
\end{bmatrix}\!,\label{eq:cha_flux}\\
x_{bdc}&=(\phi_{L},Q_{bat},Q_{\text{\tiny TL}})^\top,\nonumber \\ H_{bdc}&=\diag(L^{-1},C^{-1},C_{\text{\tiny TL}}^{-1}),\nonumber
\end{align}
where $d_{dc}\in[0,1]$ is the duty ratio used for the converter.

A simplified model for the converter can be obtained from the steady state equations of \eqref{eq:cha_flux}, which imply that  
$\overline{i_{bat}}\overline{v_{bat}}=\overline{i_{\text{\tiny TL}}}\overline{v_{\text{\tiny TL}}}$, i.e.\ the input power equals the output power. These models might as well be used for control purposes \cite{CupZM15}.

\subsubsection*{Distribution line model}
The currents and voltages of the bidirectional converter $(v_{\text{\tiny TL}},i_{\text{\tiny TL}})=(v_{bdc},i_{bdc})$ are connected to a DC voltage source with current and voltage $(v_{dc},i_{dc})$ via a distribution line which is modeled using the $\pi$-model, see \cite[Section 6.1]{Kun94}, which is a lumped parameter model and provides a good approximation for short lines. For longer transmission lines, it might be beneficial to replace this model by spatially  discretized telegrapher equations. 

The current direction is assumed to be from the voltage source towards the converter, hence, if we interconnect with the converter, we have to alternate the sign of the current. Using charge and flux variables leads to the following pH system
\begin{align}
\nonumber
\dot x_{tl}&=\begin{bmatrix}
\phantom{-}0&\phantom{-}0&1\\
\phantom{-}0&\phantom{-}0&1\\-1&-1&-R_{\text{\tiny TL}}
\end{bmatrix} H_{tl}x_{tl}+I_3\begin{bmatrix}
-i_{\text{\tiny TL}}\\ \phantom{-}i_{dc}
\\\phantom{-}v_{dc}
\end{bmatrix},\\[1ex]
y&=I_3\begin{bmatrix}
C_{\text{\tiny TL}}^{-1}&0&0\\
0&C_{dc}^{-1}&0\\0&0&L_{\text{\tiny TL}}^{-1}
\end{bmatrix}\begin{bmatrix}Q_{\text{\tiny TL}}\\Q_{cap}\\\phi_{ind}\end{bmatrix}=\begin{bmatrix}v_{\text{\tiny TL}}\\ v_{cap}\\i_{ind}\end{bmatrix} \nonumber ,\\[1ex]
x_{tl}&=(Q_{\text{\tiny TL}},Q_{cap},\phi_{ind})^\top, \label{eq:tl} \\  H_{tl}&=\diag(C_{\text{\tiny TL}}^{-1},C_{dc}^{-1},L_{\text{\tiny TL}}^{-1}).\nonumber
\end{align}

\subsubsection*{Interconnection within DC side subsystem}
We rearrange \eqref{eq:bat} so that we can replace $i_{bat}$ in \eqref{eq:cha_flux}. 
A~combination of the equations \eqref{eq:bat}, \eqref{eq:cha_flux}, \eqref{eq:tl} and using the state $x_{\text{\tiny DC}}=(\phi_L,Q_{bat},Q_{\text{\tiny TL}},Q_{cap},\phi_{ind})^T$ leads to the overall DC side system 
\begin{align}
\nonumber
\dot x_{\text{\tiny DC}}&=(J_{\text{\tiny DC}}(d_{dc})-R_{\text{\tiny DC}})H_{\text{\tiny DC}}x_{\text{\tiny DC}}
+\!\begin{bmatrix}
0&\!\!\!0&\!\!\!0\\R_{bat}^{-1}&\!\!\!0&\!\!\!0\\ 0&\!\!\!0&\!\!\!0\\0&\!\!\!1&\!\!\!0\\0&\!\!\!0&\!\!\!1
\end{bmatrix}\!\!\!\begin{bmatrix}
v_{ev}\\i_{dc}\\v_{dc}
\end{bmatrix}\!\!,\nonumber\\
\nonumber
y&=\begin{bmatrix}
0&\!\!\!R_{bat}^{-1}&\!\!\!0&\!\!\!0&\!\!\!0\\ 0&\!\!\!0&\!\!\!0&\!\!\!1&\!\!\!0\\0&\!\!\!0&\!\!\!0&\!\!\!0&\!\!\!1
\end{bmatrix}H_{\text{\tiny DC}}x_{\text{\tiny DC}}, 
\\[1ex]
&J_{\text{\tiny DC}}(d_{dc})-R_{\text{\tiny DC}}\label{full_dc} \\&=\begin{bmatrix}
0\!\!\!&\!\!\!-1&\!\!\!\tfrac{d_{dc}}{2}\!\!\!&\!\!\!0&\!\!\!0\\1&\!\!\!-R_{bat}^{-1}&\!\!\!0&\!\!\!0&\!\!\!0\\-\tfrac{d_{dc}}2&\!\!\!0&\!\!\!0&\!\!\! 0&\!\!\!\tfrac12\\
0&\!\!\!0&\!\!\!0&\!\!\!0&\!\!\!1\\0&\!\!\!0&\!\!\!-\tfrac12&\!\!\!-1&\!\!\!-R_{\text{\tiny TL}}
\end{bmatrix},\nonumber \\[1ex] H_{\text{\tiny DC}}&=\diag(L^{-1},C^{-1},2C_{\text{\tiny TL}}^{-1},C_{cap}^{-1},L_{\text{\tiny TL}}^{-1}).\nonumber
\end{align}

\subsection{AC side subsystem}
\label{sec:ac}
The AC side subsystems connect the AC grid via an LCL filter and an AC/DC converter to the DC bus. The LCL filter stands for Inductor-Capacitor-Inductor filter, corresponding to a passive filter used to reduce harmonic distortion in three-phase power systems. 
Its primary functions are to reduce harmonic content in the current and improve the power factor of the system. 

Assume that the grid voltage $v_{gabc}$ is a three-phase distortion-free sinusoidal waveform with constant amplitude and frequency, and that the PLL can always track the grid voltage angle $\theta(t)$ at the point of common coupling (PCC). In such a way,
\begin{align*}
    v_{gd}=\widehat{V}_{g},\quad
    v_{gq}=0,
\end{align*}
where $\widehat{V}_{g}$ is the amplitude of $v_{gabc}$.

The AC/DC converter has $(i_{dc},v_{dc})$ as input variables and generates $(i_s,v_s)$ as output variables. The output variables are three-phase sinusoidal waveforms as averaged modeling is used. In order to simplify future calculations, the AC waveforms are transformed to equivalent DC signals, denoted as $i_s=(i_{sd},i_{sq})$ and $v_s=(v_{sd},v_{sq}$), by applying the direct-quadrature (dq) transformation. In this manner, calculations can then be carried out on these DC quantities before performing the inverse transform to restore the original three-phase AC outcomes.

Denoting the modulation indexes by $m_q, m_d$, the converter is described in \cite[Section 5.3.3]{YazI10} by the averaged converter output 
\begin{align*}
    v_{sd}=\tfrac12m_dv_{dc},\quad v_{sq}=\tfrac12m_qv_{dc},
\end{align*}
$m_d,m_q\in[-1,1]$ and the power balance equation 
\begin{align*}
v_{dc}i_{dc}&=\tfrac{3}{2}(v_{sd}i_{sd}+v_{sq}i_{sq})\\ &=\tfrac{3}{4}(m_dv_{dc}i_{sd}+m_qv_{dc}i_{sq}).
\end{align*}
This can be rewritten in the form
\begin{align}
\label{eq:acdc}
\begin{bmatrix}
-\tfrac23i_{dc}\\v_{sd}\\v_{sq}
\end{bmatrix}=\begin{bmatrix}
0&\!\!-\tfrac12m_d&\!\!-\tfrac12m_q\\\tfrac12m_d&\!\!0&\!\!0\\\tfrac12m_q&\!\!0&\!\!0
\end{bmatrix} \!\begin{bmatrix}
v_{dc}\\i_{sd}\\i_{sq}
\end{bmatrix}.
\end{align}
Furthermore, we assume that there is a capacitor between the DC bus side and the AC/DC converter which is described by
\begin{align}
    \label{cap_conv}
i_{in}=C_{dc}\tfrac{\mathrm{d}}{\mathrm{d}t}v_{dc}+i_{dc}.
\end{align}
The LCL-filter between the AC/DC converter and the AC grid voltage input is modeled by 
\begin{align}
L_{\text{\scriptsize fs}}\tfrac{\mathrm{d}}{\mathrm{d}t}i_s&=L_{\text{\scriptsize fs}}\Omega i_s+v_s-(v_c+R_f(i_s-i_g)),\nonumber\\
C_f\tfrac{\mathrm{d}}{\mathrm{d}t}v_c&=C_f\Omega v_c+(i_s-i_g), \label{lcl_filter}\\
L_{\text{\scriptsize fg}}\tfrac{\mathrm{d}}{\mathrm{d}t}i_g&=L_{\text{\scriptsize fg}}\Omega i_g+v_c+R_f(i_s-i_g)-v_g,\nonumber
\end{align}
where $v_c$ is the voltage across the capacitor of the LCL-filter and $\Omega=\begin{bmatrix}
0&\omega\\-\omega&0
\end{bmatrix}$ for given $\omega>0$. 

A combination of \eqref{eq:acdc}, \eqref{cap_conv} and \eqref{lcl_filter} leads to the following AC side system 
\begin{align}
\nonumber
\dot x_{\text{\tiny AC}}&=(J_{\text{\tiny AC}}(m_d,m_q)-R_{\text{\tiny AC}})H_{\text{\tiny AC}}x_{\text{\tiny AC}}+\!\begin{bmatrix}
 \tfrac 23 i_{in}\\0\\0\\-v_g
\end{bmatrix}\!,\\[1ex]
x_{\text{\tiny AC}}&=(Q_{dc},\phi_s^\top, Q_c^\top, \phi_g^\top)^\top, \nonumber\\[1ex]
\label{full_ac}
y&=\begin{bmatrix}
1&0&0&0\\0&0&0&1
\end{bmatrix}H_{\text{\tiny AC}}x_{\text{\tiny AC}}=\begin{bmatrix}
v_{dc}\\i_g
\end{bmatrix},\\[1ex] \nonumber 
&J_{\text{\tiny AC}}(m_d,m_q)-R_{\text{\tiny AC}}\\ &=\begin{bmatrix}0&\!\!-\tfrac12\begin{smallbmatrix}
m_d&\!\! m_q
\end{smallbmatrix}&\!\!0&\!\!0\\\tfrac12\begin{smallbmatrix}
m_d\\ m_q
\end{smallbmatrix}&\!\!
L_{\text{\scriptsize fs}}\Omega-R_f&\!\!-I_2&\!\!R_f\\0&\!\!I_2&\!\!C_f\Omega&\!\!-I_2\\0&\!\!R_f&\!\!I_2&\!\!L_{\text{\scriptsize fg}}\Omega-R_f
\end{bmatrix}, \nonumber\\[1ex]
H_{\text{\tiny AC}}&=\diag(\tfrac 32C_{dc}^{-1},L_{\text{\scriptsize fs}}^{-1}I_2,C_f^{-1}I_2,L_{\text{\scriptsize fg}}^{-1}I_2).\nonumber
\end{align}

\subsection{Full model of the uncontrolled charging station}\label{sec:fullmodel}
To build the overall system of the charging station we have to interconnect the subsystems \eqref{full_dc} and \eqref{full_ac}. First observe that the currents in both subsystems have opposite directions. Hence, we have to use $i_{in}=-i_{ind}$ and $v_{in}=v_{dc}$. This leads to an overall (averaged) pH system with
\begin{align}
\nonumber
&~~~~J(d_{dc},m_d,m_q)-R\\&=\begin{bmatrix}
J_{\text{\tiny DC}}(d_{dc})-R_{\text{\tiny DC}}&\tfrac23 e_{5}e_{1}^\top
\\ \nonumber
-\tfrac23e_{1}e_{5}^\top&J_{\text{\tiny AC}}(m_d,m_q)-R_{\text{\tiny AC}}
\end{bmatrix},\\ &Bu=\begin{bmatrix}0\\R_{bat}^{-1}v_{ev}\\0_6\\-v_g\end{bmatrix},\quad  y=\begin{bmatrix}
R_{bat}^{-1}v_{ev}-i_{ev}\\i_{g}
\end{bmatrix},\label{eq:full_charge}\\
&x=(\phi_L,Q_{bat},Q_{bdc},\phi_{ind},Q_{dc},\phi_{s}^\top,Q_{c}^\top,\phi_{g}^\top)^\top,\nonumber
\end{align}
where $v_{dc}$ appears as an additional state variable, $e_5\in\R^5$ and $e_1\in\R^7$ are canonical unit vectors and $\mathcal{H}=\mathcal{H}_{\text{\tiny DC}}+\mathcal{H}_{\text{\tiny AC}}$. For the charging station we have three parameters $\theta_1=d_{dc}$, $\theta_2=m_d$, $\theta_3=m_q$ and the matrices in \eqref{def_sys} are given by 
\begin{align*}
J_1&=\tfrac12(e_1e_3^\top-e_3e_1^\top),\ J_2=\tfrac12(e_6e_5^\top-e_5e_6^\top),\\ J_3&=\tfrac12(e_7e_5^\top-e_5e_7^\top),
\end{align*}
where again $e_1,e_3,e_5,e_6,e_7\in\R^9$ are canonical unit vectors.

\section{Numerical simulation of an EVCS and comparison of controllers}
\label{sec:comp}

In this section, we compare the pH approaches proposed in  Theorem~\ref{prop:control} for controlling the duty ratios and modulation indexes of the averaged pH model of a particular EVCS. Furthermore, the pH controllers are verified for a switched system model of the EVCS in Simulink. Here we use the following steady state duty ratios which were obtained from the Simulink model of the charging station 
\begin{align*}
\overline{\theta_1}&=\overline{d_{dc}}=\tfrac{5}{9}, &&\overline{\theta_2}=\overline{m_d}=0.726,\\ \overline{\theta_3}&=\overline{m_q}=-0.018, &&
\end{align*}
and, hence, it can be verified that  $J(\overline{\theta})-R$ given by \eqref{eq:full_charge} is Hurwitz. Thus, the steady state given by \eqref{eq:ss_inv} is asymptotically stable, and the control approaches described in Theorem~\ref{prop:control} lead to asymptotically stable pH systems.

To compare the controllers described in Theorem~\ref{prop:control}, we will consider a perturbation of the zero initial condition of the shifted system, which means that we perturb the steady state solution and investigate the settling time of the controllers. 

For the derivation of the controllers, we compute 
\begin{align*}
-\mathcal{D}(\overline{x})&=\begin{bmatrix}-2\overline{v_{bdc}}&0&0\\
0&0&0\\
\overline{i_L}&0&0\\
0&0&0\\
0&\tfrac12\overline{i_{sd}}&\tfrac12\overline{i_{sq}}\\
0&-\tfrac34\overline{v_{dc}}&0\\
0&0&-\tfrac34\overline{v_{dc}}\\
0&0&0\\0&0&0
\end{bmatrix}
,\\ \nabla\mathcal{S}(x)&=[i_L-\overline{i_L}~~~~~2(v_{bat}-\overline{v_{bat}})~~v_{bdc}-\overline{v_{bdc}}\\&~~~~i_{ind}-\overline{i_{ind}}~~\tfrac32(v_{dc}-\overline{v_{dc}})~~~~i_{s}-\overline{i_s}\\&~~~~i_{c}-\overline{i_c}~~~i_{g}-\overline{i_g}]^\top,\\
K_d&=\gamma I_3 \quad \text{for some $\gamma>0$}.
\end{align*}

Applying the proportional pH controller introduced in \eqref{eq:ext_ode} to the charging station model results in the following equations 
\begin{align}
\dot d_{dc}&=2(v_{bdc}\overline{i_L}-\overline{v_{bdc}}i_L)-\gamma(d_{dc}-\overline{d_{dc}}),\nonumber\\
\dot m_d&=\tfrac34(\overline{i_{sd}}v_{dc}-\overline{v_{dc}}i_{sd})-\gamma(m_d-\overline{m_d}),\label{eq:phODE}\\
\dot m_q&=\tfrac34(\overline{i_{sq}}v_{dc}-\overline{v_{dc}}i_{sq})-\gamma(m_q-\overline{m_q}), \nonumber
\end{align}
where $\gamma>0$ is a control gain factor, that is chosen sufficiently large in order to compensate the size of the two summands written in parentheses in \eqref{eq:phODE}. Here, the first summand is a difference of powers which typically attains large values of approximately $10^3$ and $10^4$ for the considered power system and the second summand with prefactor $\gamma$ is a difference of the rather small duty ratios or modulation indexes, whose absolute difference is trivially bounded from above by $2$.

Furthermore, the (differential-algebraic) pH controller introduced in \eqref{eq:ext_dae} is given by \eqref{eq:phODE} after replacing the left-hand sides by $0$.

Moreover, we consider the pH PI controller which is given by equations
\begin{align}
\nonumber
0&=\delta\cdot2(v_{bdc}\overline{i_L}-\overline{v_{cap}}i_L)-\gamma(d_{dc}-\overline{d_{dc}})+\xi_1,\\ \nonumber
0&=\delta\cdot\tfrac34(\overline{i_{sd}}v_{dc}-\overline{v_{dc}}i_{sd})-\gamma(m_d-\overline{m_d})+\xi_2,\\ \nonumber
0&=\delta\cdot\tfrac34(\overline{i_{sq}}v_{dc}-\overline{v_{dc}}i_{sq})-\gamma(m_q-\overline{m_q})+\xi_3,\\ \label{eq:phPI}
\dot \xi_1&=\gamma^{-1}2(v_{bdc}\overline{i_L}-\overline{v_{bdc}}i_L),\\\nonumber 
\dot \xi_2&=\gamma^{-1}\tfrac34(\overline{i_{sd}}v_{dc}-\overline{v_{dc}}i_{sd}),\\ 
\dot \xi_3&=\gamma^{-1}\tfrac34(\overline{i_{sq}}v_{dc}-\overline{v_{dc}}i_{sq}),\nonumber
\end{align}
where a suitable choice of $\delta>0$ improves the performance of the controller, in the sense that the simulated time of convergence towards the given steady state is short, which is a standard measure to evaluate the performance of a stabilizing controller. Another measure that is often used is the size of the overshoots after set-point changes, which we can see, for instance, in the Figures~\ref{fig:averaged_simulation}, right at the beginning of the simulated time-period.

\subsection{Comparison of the pH controllers for the averaged model}

First, the different pH controllers are applied and compared to the averaged model of an EVCS. For the simulation, we assumed for the following states an initial deviation from the prescribed set points 
\begin{align*}
|i_{L}(0)-\overline{i_L}|&=50\,\text{A},& |v_{bat}(0)-\overline{v_{bat}}|&=50\,\text{V},\\ |v_{dc}(0)-\overline{v_{dc}}|&= 50\,\text{V},& |i_{sq}(0)-\overline{i_{sq}}|&= 75\,\text{A},
\end{align*}
whereas all other states are initialized with the corresponding value of the set-point.

\begin{figure}
    \includegraphics[scale=0.55]{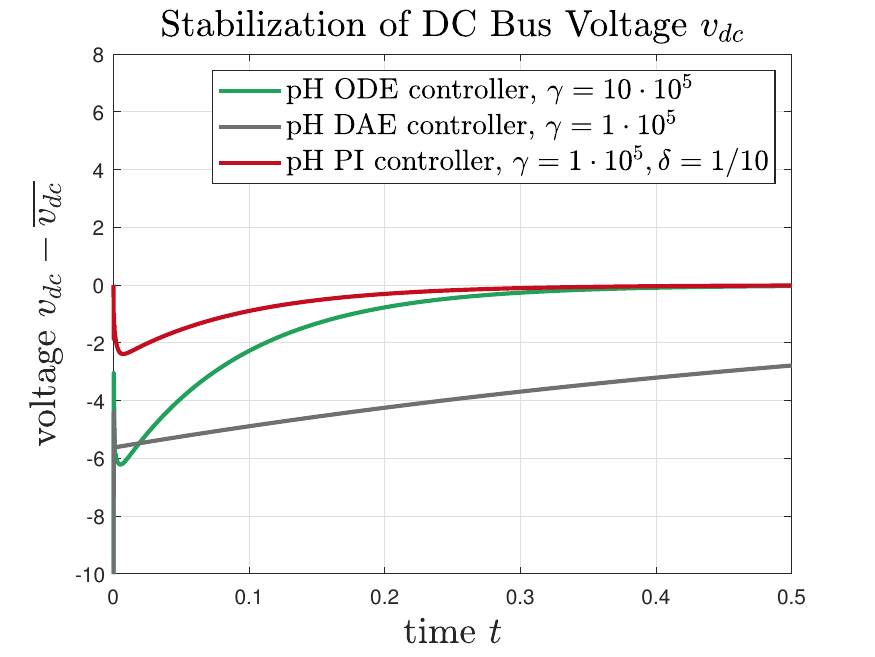} \includegraphics[scale=0.55]{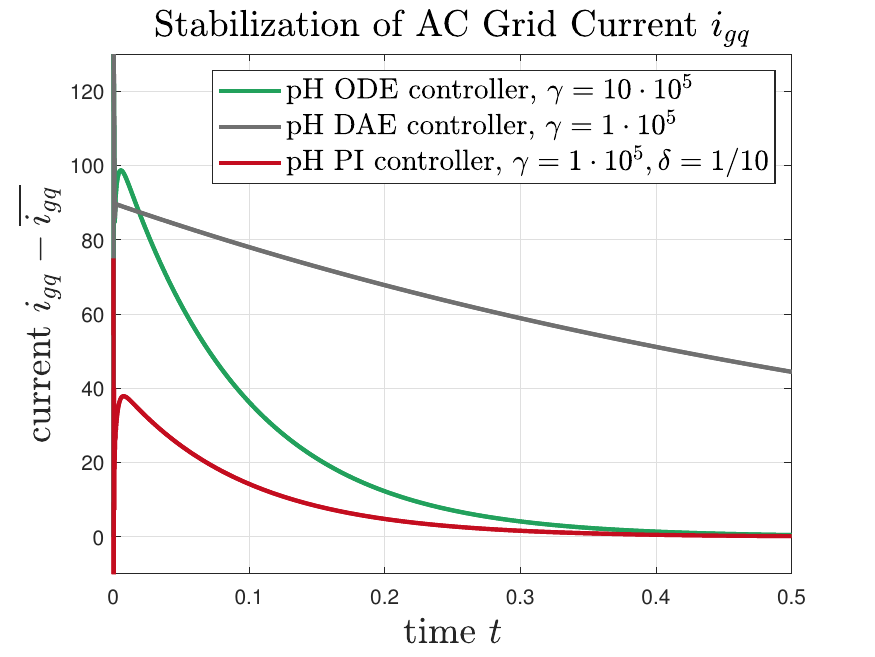}
    \caption{Simulation results of different pH controllers applied to the averaged model~\eqref{eq:full_charge} of an EVCS connected to an AC grid. The plots show the trajectories of the pH controlled closed-loop system, where the DC bus voltage $v_{dc}$ and the $q$-component of the AC grid current $i_{gq}$ were selected.}
    \label{fig:averaged_simulation}
\end{figure}

The methodology of choosing the control gains $\gamma,\delta>0$ is as follows: we choose a sufficiently large value of $\gamma=10^5$ based on the order of magnitude of the coefficients in \eqref{eq:phODE} and then the additional parameter $\delta$ for the pH PI controller is chosen in a trial-and-error process by starting with some initial value $\delta_0>0$ and then, based on the outcome of the simulation results, the current value of $\delta$ is either increased or decreased. As a performance measure for this decision, we use the speed of convergence towards the steady state. Another performance measure that is typically used in control  design is to bound the overshoot, which can also be observed in Figure~\ref{fig:averaged_simulation}. The simulation results show that the pH PI controller outperforms the pH ODE and the pH DAE controllers regarding the convergence time and the overshooting behavior.

Furthermore, when increasing the control gain $\gamma$ it can be observed that the convergence time of the pH ODE and pH DAE controller increases. To illustrate this in Figure~\ref{fig:averaged_simulation}, we chose a larger control gain $\gamma>0$ for the pH ODE controller compared to the pH DAE controller. For the same control gain, the performance of the pH ODE and the pH DAE controller is similar. Nevertheless, merely increasing the control gain $\gamma$ in the pH ODE or pH DAE controllers cannot reduce the overshooting behavior of the $q$-component of the AC-grid current $i_{gq}$ and therefore using the pH PI controller is the better option.

\subsection{Validation of the pH controllers}
We validate the performance of the designed pH controllers by performing time-domain simulations of the switching model of the EVCS system, as described in Section~\ref{sec:charging}.

We compare the pH controllers of the DC/DC converter and DC/AC inverter with their corresponding proportional-integral (PI) type controllers of the DC/DC converter and DC/AC inverter, respectively. The cascaded PI controller, as described in \cite{AranStru12}, is detailed in Section~\ref{sec:app_cascaded_PI} of the appendix. In principle, this controller adopts a cascaded configuration in the DC/AC inverter, where the outer voltage controller controls the DC bus voltage $v_{dc}$ and generates the reference of $i_{sd}$, and the inner current controller controls the current $i_{sd}$ and $i_{sq}$.

For the pH PI controller that is given by~\eqref{eq:phPI}, we use the control gains $\delta=1.5$ and $\gamma=3\cdot 10^{6}$. Their values are derived from the simulations using the averaged model and then they are adapted in the above described trial-and-error process to increase the performance of the controllers.

\begin{figure}[!htb]
        \begin{minipage}[t]{.9\linewidth}
        \subfigure[]{   
        \includegraphics[scale=0.11]{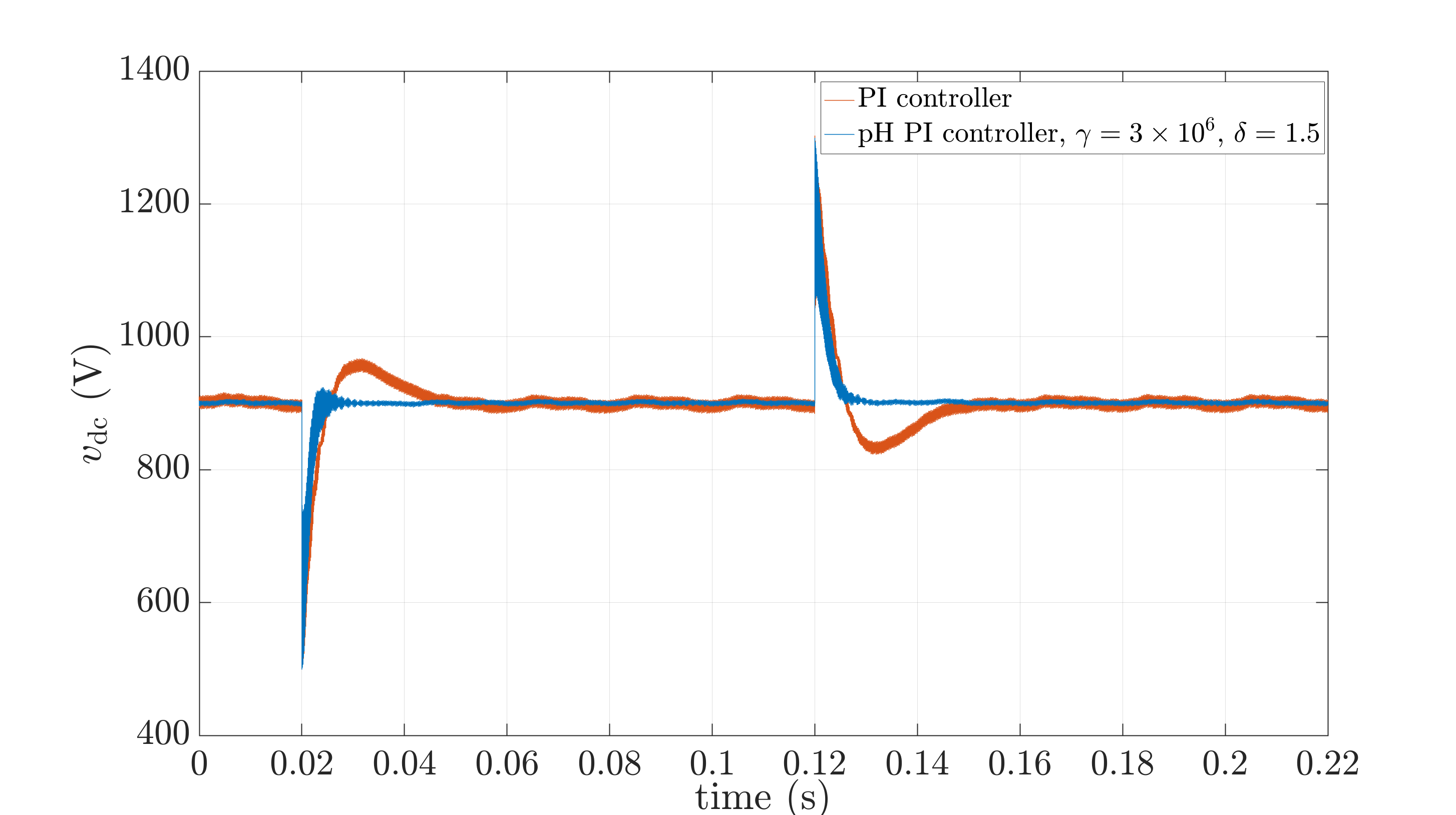}}    
         \label{vdc}
        \end{minipage}
        \begin{minipage}[t]{0.9\linewidth}
        \subfigure[]{ 
        \includegraphics[scale=0.11]{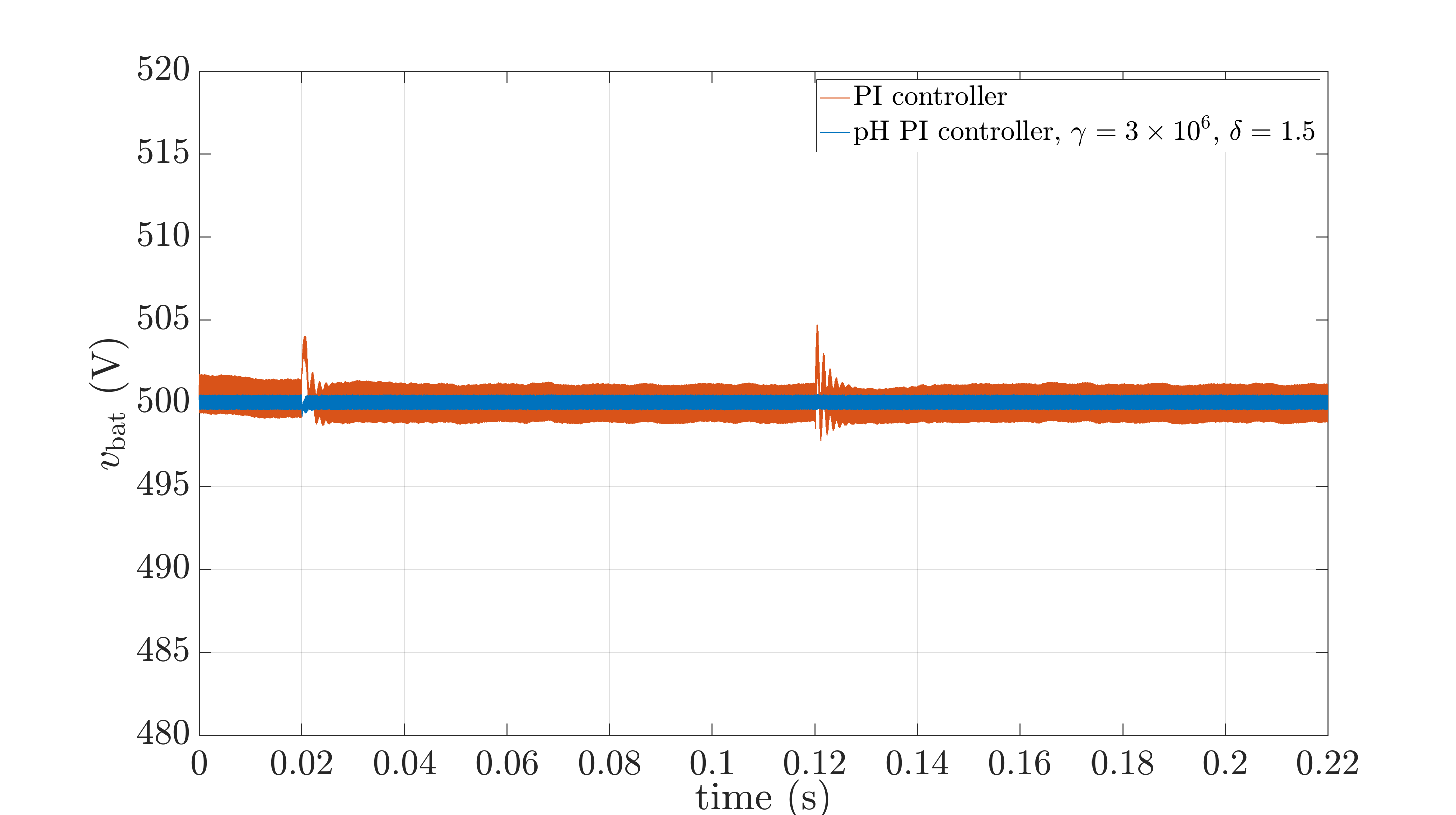} 
        }
        \label{vbat}
        \end{minipage}
        \begin{minipage}[t]{0.9\linewidth}
        \subfigure[]{ 
        \includegraphics[scale=0.11]{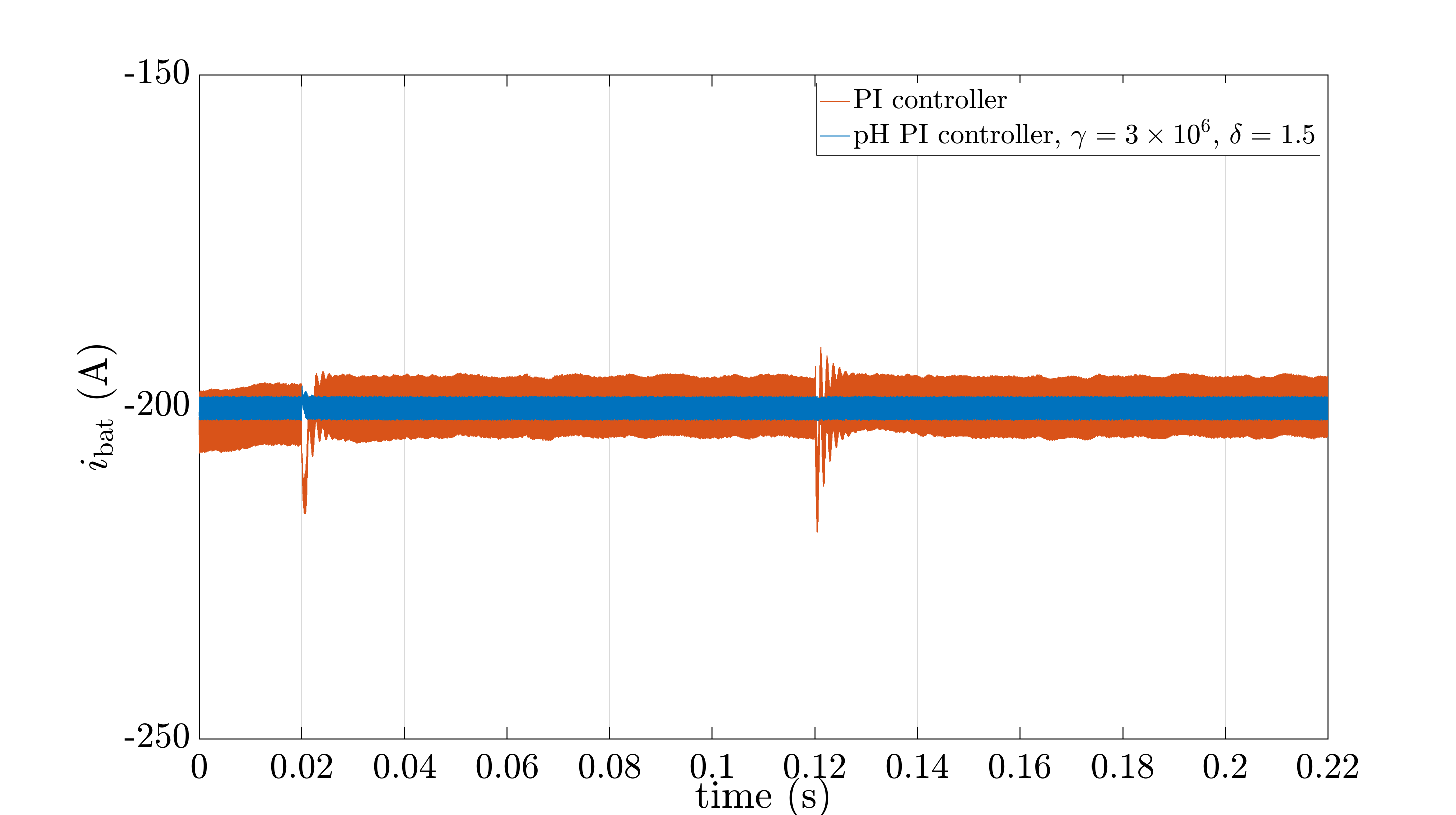} 
        }
        \label{ibat}
        \end{minipage}
        \caption{Simulation results of the EV charging station switched model with the pH PI and the cascaded PI controller: (a) DC bus voltage; (b) output voltage of the EV battery; (c) output current of the EV battery.  
        }    
        \label{simulation}    
    \end{figure}

Furthermore, compared to the charging station model considered in \cite{AranStru12}, an additional power line with $L_{\text{\tiny TL}}=1.1\cdot 10^{-6}$\,H and $C_{\text{\tiny TL}}=4\cdot 10^{-4}$\,F and $R_{\text{\tiny TL}}=10^{-3}$ is introduced in the generic EVCS model that is described in Section~\ref{sec:charging}. 
The references of $v_{bat}$, $v_{dc}$, and $i_{sq}$ used in the simulations are
\begin{align*}
v_{bat}^{\text{ref}}&=\overline{v_{bat}}=500\,\text{V},\quad v_{dc}^{\text{ref}}=\overline{v_{dc}}=900\,\text{V},\\ i_{bat}^{\text{ref}}&=-200\,\text{A}.
\end{align*}
In the simulation scenario, all variables within the charging station remain at their respective steady-states from 0 s to 0.02 s. Subsequently, during the time interval from 0.02 s to 0.12 s, disturbances in the form of changes to $v_{dc}$ is introduced by setting 
\begin{align*}
  v_{dc}&=(\overline{v_{dc}}-500)\,\text{V}.  
\end{align*}
Moreover, at $0.12$ s, the disturbances are removed. The dynamics of $v_{dc}$ is shown in Figure~\ref{simulation}. 

The simulation results in Figure~\ref{simulation} show that the considered pH controllers can handle large disturbances in the DC bus voltage better than the cascaded PI controller.

Furthermore, we observed that the line parameters of the system might behave sensitive to perturbations. It was noticed that after increasing the line parameters $L_{\text{\tiny TL}}$ and $C_{\text{\tiny TL}}$ significantly, by one order of magnitude the set-point changes might cause some considerable resonance in the system. Although both control approaches are still able to track the set-points, this change of parameters has a negative impact on the performance of both, the pH PI and the cascaded PI controller.

\section{Conclusion and future work}
We propose averaged pH systems as a modeling framework that allows to include averaged models of switched subsystems, such as power converters. Using the pH modeling framework, we derive several types of proportional and PI controllers for averaged pH systems in general and for EV charging stations in particular. The design of the controllers is only based on the knowledge of the system equations of the averaged models, making them easily adaptable to different models of EV charging stations or larger systems of interconnected charging stations by adapting the two scalar parameters of the pH PI controller. The designed controllers were successfully validated in simulations using the switched model of a particular EV charging stations, and the advantages compared to a standard cascaded PI controller are discussed.

In the future, we plan to consider optimized pH control designs, as there is still a considerable amount of freedom in the choice of the control design matrices which has not been explored yet.

\section*{Acknowledgment}
HG, KS, VM, and XZ gratefully acknowledge the funding within the SPP1984 ``Hybrid and multimodal energy systems'' by the Deutsche Forschungsgemeinschaft (DFG) project number 361092219. Furthermore, the authors would like to thank the anonymous referees for their careful reading and suggestions to improve the overall quality of the manuscript.




\bibliographystyle{IEEEtran}
\bibliography{phcharging}

\begin{appendix}
    
\section{Proof of Theorem \ref{prop:control}}
    \label{sec:app_control}
To prove (a), observe that 
\eqref{eq:minusbar} with $\theta=\overline{\theta}$ is equivalent to
$\dot x=D(\overline{\theta})\nabla\mathcal{S}(x).$ 
To prove that the system in (b) is pH, we consider the following decomposition of the coefficients in \eqref{eq:ext_ode}:
\begin{align}
\label{dec_sym_skew}
\begin{bmatrix}D(\theta)&\!\!\!\!\!\mathcal{D}(\overline{x}) \\ -\mathcal{D}(\overline{x})^\top & \!\!\!\!\!-K_{\theta} \end{bmatrix}\!=\!\begin{bmatrix}J(\theta)&\!\!\!\!\!\mathcal{D}(\overline{x}) \\ -\mathcal{D}(\overline{x})^\top  & \!\!\!0\end{bmatrix}-\begin{bmatrix} R(\theta)&\!\!\!\!\!0\\0&\!\!\!\!\!K_{\theta}\end{bmatrix}.
\end{align}
By assumption that $\theta_i\geq0$, the symmetric part is positive semi-definite and therefore the system is pH. Next, we show that the system is locally asymptotically stable around the equilibrium $(\overline{x},\overline{\theta},\overline{u})$ by computing the Jacobian $\mathrm{D}$ of the system in $(\overline{x},\overline{\theta})$ which is given by 
\begin{align*}
&\mathrm{D}\left(\begin{bmatrix} D(\theta) & \mathcal{D}(\overline{x})\\  -\mathcal{D}(\overline{x})^\top&K_{\theta}\end{bmatrix} \nabla\mathcal{S}^e(x,\theta)\right)|_{(x,\theta)=(\overline{x},\overline{\theta})}\\&=\begin{bmatrix}D(\overline{\theta})&\mathcal{D}(\overline{x}) \\ -\mathcal{D}(\overline{x})^\top  & -K_{\theta} \end{bmatrix}.
\end{align*}
This matrix is Hurwitz, since $K_{\theta}$ is positive definite and its Schur complement is given by the Hurwitz matrix  $D(\overline{\theta})-\mathcal{D}(\overline{x})K_\theta^{-1}\mathcal{D}(\overline{x})^\top$. Hence, the system is locally asymptotically stable. 

We continue with the proof of (c). To obtain an equivalent standard system, we solve the algebraic equation for $\theta-\overline{\theta}$ as follows:
\begin{align*}
\theta-\overline{\theta}&=\delta K_{\theta}^{-1}K_H\nabla\mathcal{S}(x),\\ \theta_i-\overline{\theta}_i&=\delta e_i^\top K_{\theta}^{-1}K_H\nabla\mathcal{S}(x),
\end{align*}
where $e_i\in\mathbb{R}^k$ denotes the canonical unit vector which has an entry $1$ in the $i$th row and all other entries are zero. 
Hence, the first line in \eqref{eq:ext_dae} can be reformulated as follows
\begin{align}
\dot x&=(D-\sum_{j=1}^k\theta_jD_j)\nabla\mathcal{S}(x)+\mathcal{D}(\overline{x})(\theta-\overline{\theta})\nonumber \\&~~~+ B(u-\overline{u})\nonumber\\
&=\Big(D(\overline{\theta})-\delta\sum_{j=1}^k(e_j^\top K_{\theta}^{-1}K_H\nabla\mathcal{S}(x))D_j\nonumber \\&~~~+\mathcal{D}(\overline{x})K_{\theta}^{-1}K_H\Big)\nabla\mathcal{S}(x)+ B(u-\overline{u}) \label{def_D}\\
&=F(x)\nabla\mathcal{S}(x)+ B(u-\overline{u}),
\nonumber
\end{align}
where we abbreviate the matrix function in the left bracket in \eqref{def_D} by $F$. If $K_H=-\mathcal{D}(\overline{x})^\top$ then $F$ can be decomposed as 
\begin{align}
F(x)&=\hat J(x)-\hat R(x),\quad \hat J(x)=-\hat J(x)^\top \nonumber \\\nonumber \hat R(x)&=R(\theta)+\delta\sum_{j=1}^k(e_j^\top K_{\theta}^{-1}K_H\nabla\mathcal{S}(x))R_j\\&+\mathcal{D}(\overline{x})K_{\theta}^{-1}\mathcal{D}(\overline{x})^\top\geq 0, \label{eq:D_dec}\\ \hat J(x)&=J(\overline{\theta})-\delta\sum_{j=1}^k(e_j^\top K_{\theta}^{-1}K_H\nabla\mathcal{S}(x))J_j .\nonumber
\end{align}
Hence, the following system is a standard pH system of the form \eqref{pH_DAE}:
\begin{align}
\dot x=F(x)\nabla \mathcal{S}(x)+ B(u-\overline{u}),\ 
y=  B^\top \nabla \mathcal{S}(x).\label{eq:PI_equiv}
\end{align}
The local asymptotic stability can be proven by computing the Jacobian of $F(x)\nabla\mathcal{S}(x)$ in $x=\overline{x}$ which is equal to  $D(\overline{\theta})-\mathcal{D}(\overline{x})K_{\theta}^{-1}\mathcal{D}(\overline{x})^\top$ and thus a Hurwitz matrix by Proposition~\ref{lem:Hurwitz}. 

In an analogous way, and using the abbreviation
\begin{align*}
\hat F(x,\xi)&:=D(\overline{\theta})+\mathcal{D}(\overline{x})K_{\theta}^{-1}K_H
\\&-\delta\sum_{j=1}^k(e_j^\top( K_{\theta}^{-1}K_H\nabla\mathcal{S}(x)+K_{\theta}^{-1}\xi))D_j,
\end{align*}
system \eqref{eq:PI} can be rewritten as the standard system
\begin{align}
\label{eq:PI_ode}
\begin{split}
\begin{bmatrix}
\dot x\\ \dot \xi
\end{bmatrix}&\!=\!\begin{bmatrix} \hat F(x,\xi)& \!\!\!\mathcal{D}(\overline{x})K_{\theta}^{-1}\\K_{\xi}&\!\!\!0\end{bmatrix}\!\!\begin{bmatrix}\nabla \mathcal{S}(x)\\\xi\end{bmatrix}\!\!+\!\!\begin{bmatrix}B\\0\end{bmatrix}(u-\overline{u}),\\
y&\!=\!B^\top \begin{bmatrix}
\nabla \mathcal{S}(x)\\ \xi \end{bmatrix}. 
\end{split}
\end{align}
Clearly, the matrix $\hat F(x,\xi)$ is Hurwitz in a neighborhood of $(\overline{x},0)$ and has positive definite symmetric part. Hence,  
for the choice $K_{\xi}=-(\mathcal{D}(\overline{x})K_{\theta}^{-1})^\top$, the system \eqref{eq:PI_ode} is locally asymptotically stable and pH. 

    \section{Proof of Proposition \ref{lem:Hurwitz}} \label{sec:app_hurwitz}
Assume that there exists an eigenvalue  $\lambda\in\mathbb{C}$ of $A$ with eigenvector $x\in\mathbb{C}^n\setminus\{0\}$ and its entry-wise complex conjugate $\overline{x}$. Then
\begin{align*}
0\geq -2\overline{x}^\top(R+\mathcal{D}K_{\theta}^{-1}\mathcal{D}^\top)x&=\overline{x}^\top Ax+\overline{x}^\top A^\top x\\&=2\re\lambda\|x\|^2.
\end{align*}
Hence, $\re\lambda\leq 0$ holds. If $\re\lambda=0$ then the positive semi-definiteness of $R$ and $\mathcal{D}K_{\theta}^{-1}\mathcal{D}^\top$ implies that 
\begin{align*}
x\in\ker(R+\mathcal{D}K_{\theta}^{-1}\mathcal{D}^\top)&=\ker R\cap\ker \mathcal{D}K_{\theta}^{-1}\mathcal{D}^\top\\&=\ker R\cap\ker K_{\theta}^{-1/2}\mathcal{D}^\top\\&=\ker R\cap\ker \mathcal{D}^\top\\&=\{0\}
\end{align*}
holds, where we used in the last equality that $\rk\begin{smallbmatrix}
R\\\mathcal{D}^\top
\end{smallbmatrix}=\rk[R,\mathcal{D}]=n$. As a consequence, $x=0$ which implies that any eigenvalue $\lambda$ of $A$ fulfills $\re\lambda<0$ and therefore $A$ is Hurwitz.

\section{Cascaded PI controller}
\label{sec:app_cascaded_PI}
In the simulations, we use the following cascaded PI controller described in \cite{AranStru12} 
\begin{align}
\dot\xi_4&=\overline{v_{bat}}-v_{bat},\\ 2d_{dc}&=k_p^{EV}(\overline{v_{bat}}-v_{bat})+k_i^{EV}\xi_4,\nonumber \\
\dot \xi_1&=v_{dc}-\overline{v_{dc}},\nonumber \\
\dot\xi_2&=k_p^{vc}(\overline{v_{dc}}-v_{dc})+k_i^{vc}\xi_1-i_{sd},\nonumber \\
m_d&=\frac{k_p^{cc}\dot\xi_2+k_i^{cc}\xi_2-L_{fi}\omega i_{sq}+v_{cd}}{\tfrac12\overline{v_{dc}}}, \label{eq:cascaded_PI}\\
\dot \xi_3&=\overline{i_{sq}}-i_{sq},\nonumber \\
m_q&=\frac{k_p^{cc}(\overline{i_{sq}}-i_{sq})+k_i^{cc}\xi_3 +L_{fi}\omega i_{sd}+v_{cq}}{\tfrac12\overline{v_{dc}}}, \nonumber
\end{align}
with parameters
\begin{align*}
k_p^{vc}&=0.64,\quad k_i^{vc}=70,\quad k_p^{cc}=\tfrac16,\quad k_i^{cc}=\tfrac43,\\
k_p^{EV}&=1.5,\quad k_i^{EV}=10.
\end{align*}
These parameters were tuned based on the linearized model and an analysis of the eigenvalue behavior under parameter variations.

\end{appendix}

\end{document}